\documentclass{amsart}
\usepackage{amssymb,stmaryrd}
\usepackage{amsfonts}
\usepackage{amstext}
\usepackage{algorithmic}
\usepackage{algorithm}
\usepackage{graphicx}
\usepackage{epstopdf}
\usepackage[all]{xy}

\numberwithin{equation}{section}
\newtheorem{theorem}{Theorem}[section]
\newtheorem{lemma}[theorem]{Lemma}
\newtheorem{proposition}[theorem]{Proposition}
\newtheorem{corollary}[theorem]{Corollary}

\theoremstyle{definition}

\newtheorem{example}[theorem]{Example}

\newtheorem{conjecture}[theorem]{Conjecture}
\theoremstyle{remark}

\newcommand{\Z}{{\mathbb{Z}}}

\newcommand{\F}{{\mathbb{F}}}
\newcommand{\N}{{\mathbb{N}}}

\newcommand{\<}{{\langle}}

\renewcommand{\>}{{\rangle}}
\newcommand{\Tr}{{\rm Trace}}

\newcommand{\CF}{{\mathcal{F}}}

\newcommand{\CI}{{\mathcal{I}}}
\newcommand{\CJ}{{\mathcal{J}}}

\renewcommand{\ker}{{\rm{ker}}}

\newcommand{\tens}{\otimes}

\newcommand{\id}{{\rm id}}

\renewcommand{\o}{{}_{(1)}}
\renewcommand{\t}{{}_{(2)}}
\renewcommand{\th}{{}_{(3)}}
\newcommand{\extd}{{\rm d}}
\newcommand{\del}{{\partial}}
\newcommand{\eps}{\epsilon}

\newcommand{\uo}{{{}^1}}
\newcommand{\ut}{{{}^2}}

\begin{document}

\title{Finite noncommutative geometries related to $\F_p[x]$}


\author{M.E. Bassett \& S. Majid}
\address{Queen Mary, University of London\\
School of Mathematics, Mile End Rd, London E1 4NS, UK}

\email{ s.majid@qmul.ac.uk}


\begin{abstract} It is known that irreducible noncommutative differential structures over $\F_p[x]$ are classified by irreducible monics $m$. We show that the cohomology $H_{\rm dR}^0(\F_p[x]; m)=\F_p[g_d]$ if and only if $\Tr(m)\ne 0$, where $g_d=x^{p^d}-x$ and $d$ is the degree of $m$. This implies that there are ${p-1\over pd}\sum_{k|d, p\nmid k}\mu_M(k)p^{d\over k}$ such noncommutative differential structures ($\mu_M$ the M\"obius function). Motivated by killing this zero'th cohomology, we consider the directed system of finite-dimensional Hopf algebras $A_d=\F_p[x]/(g_d)$ as well as their inherited bicovariant differential calculi $\Omega(A_d;m)$. We show  that $A_d=C_d\tens_\chi A_1$ a cocycle extension where $C_d=A_d^\psi$ is the subalgebra of elements fixed under $\psi(x)=x+1$. We also have a Frobenius-fixed subalgebra $B_d$ of dimension $\frac{1}{d} \sum_{k | d} \phi(k) p^\frac{d}{k}$ ($\phi$ the Euler totient function), generalising Boolean algebras when $p=2$. As special cases, $A_1\cong \F_p(\Z/p\Z)$, the algebra of functions on the finite group $\Z/p\Z$, and we show dually that $\F_p\Z/p\Z\cong\F_p[L]/(L^p)$ for a `Lie algebra' generator $L$ with $e^L$ group-like, using a truncated exponential. By contrast,  $A_2$ over $\F_2$ is a cocycle modification of $\F_2((\Z/2\Z)^2)$ and is a 1-dimensional extension of the Boolean algebra on 3 elements. In both cases we compute the Fourier theory, the invariant metrics and the Levi-Civita connections within bimodule noncommutative geometry.   \end{abstract}
\maketitle 
\section{Introduction} 

This article is motivated by a fundamental issue in characteristic $p>0$ geometry visible even for polynomials $\F_p[x]$ in one variable over the finite field of order $p$, namely the failure of differential calculus to provide an effective tool. Specifically, on any connected manifold the only functions killed by the exterior derivative are the constant functions, i.e. the zeroth de Rham cohomology $H^0_{\rm dR}$ is spanned by $1$. By contrast, the classical differential calculus on $\F_p[x]$ has a large kernel for $\extd$, namely all polynomials in $x^p$. One approach is to quoient out this kernel to give the Hopf algebra $\F_p[x]/(x^p)$ and one will then have that the inherited calculus on this is now connected. On the other hand, this Hopf algebra is too small to serve as an approximation of $\F_p[x]$. We ask if we can do rather better by looking not at the usual differential calculus but a noncommutative one. 

We first recall the usual K\"ahler differential for commutative algebras $A$. This consists of a left module $\Omega^1$ and a map $\extd$ to it universal with the derivation property $\extd(ab)=a.\extd b+ b.\extd a$ for all $a,b\in A$, and can be built explicitly on $\CJ/\CJ^2$ where $\CJ=\ker(\cdot: A\tens A\to A)$ and $\extd a=1\tens a-a\tens 1$, see for example \cite{Kah}. For $A=k[x]$, this recovers its usual differential calculus. These ideas adapt to the case of noncommutative $A$, the key being to replace $\Omega^1$ by a bimodule and the derivation rule by $\extd(ab)=a.\extd b+(\extd a).b$. In this case $\CJ$ itself with the same $\extd$ as before provides the universal calculus and now makes sense for noncommutative $A$. If $A$ happens to be commutative then the universal calculus is now much bigger than before and generally has many interesting quotients beyond the K\"ahler one. These will typically have differentials noncommuting with elements of $A$ even though $A$ itself may be commutative. We will also need higher differential forms $\Omega$ forming a graded algebra with $\extd$ extended by a similarly two-sided  (now graded) Leibniz rule and obeying $\extd^2=0$, i.e. a differential graded algebra (or DGA). Such a notion features in most approaches to noncommutative geometry, including \cite{Con} (athough there not as a  starting point). We refer to the cohomology of the complex $(\Omega,\extd)$ as the `noncommutative de Rham cohomology' $H_{\rm dR}(A)$. 

In this paper we will be interested in the case where $A$ is a Hopf algebra, which we think of as if functions on a group (though it does not have to be, even when $A$ is commutative). Then `group translation' is expressed by the Hopf algebra coproduct $\Delta$ viewed as a coaction from the left or the right. A differential calculus is covariant if one or both of these coactions extend to $\Omega^1$. This situation has been extensively studied starting with \cite{Wor} and it is known that  $\Omega^1$ in the left-covariant case is free as a left $A$-module, having a form isomorphic to $A\tens\Lambda^1$ where $\Lambda^1$ is the space of left-invariant differential forms. In the bicovariant case, there is a canonical extension of  $\Lambda^1$ to a `braided exterior algebra' $\Lambda$ and hence of $\Omega^1$ to a DGA $\Omega$ constructed as a  Radford biproduct or `cobosonisation' $A\ltimes\Lambda$. Hence we need only focus on the choice of $\Omega^1$. We refer to \cite{Ma:prim,Ma:ltcc,Ma:hod} for more details and an introduction to what is now a large literature. 

In particular, we can consider $k[x]$ as  a Hopf algebra with $x$ `primitive' in the sense $\Delta x=x\tens 1+1\tens x$ (i.e. the additive group structure of the affine line). Then \cite{Ma:fie,Ma:prim}  irreducible bicovariant calculi on $k[x]$ (in the sense of having no proper quotients) correspond to monic irreducible $m\in k[x]$ and take the following form. First, define the associated field extension $K=k[\mu]/(m(\mu))$ and set $\Omega^1(k[x]; m)=K[x]$ as a left $k[x]$-module in the obvious way. The right module structure and differential are then
 \[  v\cdot f(x)=f(x+\mu)v,\quad \extd f=(f(x+\mu)-f(x))\mu^{-1},\quad\forall f\in k[x],\ v\in K\]
where expressions on the right are written in terms of the algebra $K[x]$ and $\mu\in K$. If $m$ has degree 1 then $K$ can be identified with $k$ and the relation just sets $\mu$ to be an element of $k$, including the case $m=x$ or $\mu=0$ as the classical commutatative calculus (the formula for the differential still makes sense in spite of appearances). Here $K$ as a vector space over $k$ is the space of left-invariant 1-forms and the canonical  $\Omega(k[x];m)$ is a free module over its usual exterior algebra. 

Our main result (Theorem~\ref{H0d}) is that for a calculus on $\F_p[x]$ defined by monic irreducible $m(x)$ of degree $d\ge 1$, the zeroth cohomology $H^0_{\rm dR}(\F_p[x];m)$ with respect to the calculus defined by $m$ consists precisely of polynomials in  $g_d=x^{p^d}-x$  if and only if  the number-theoretic `trace' of $m$ is nonzero. We say that such $m$ are `regular' and the result implies that the cohomology is independent of $m$ in this case. We also have a conjecture for the cohomology in the non-regular case which our current methods do not prove but which we also believe to be true (Conjecture~\ref{H0conj}).  The higher $H^i_{\rm dR}(\F_p[x];m)$ for $i>0$ and the canonical $\Omega(k[x];m)$ remain mysterious even to the point of  a conjecture, but are known to be nontrivial in degree 1 provided $m\ne x$, and are expected on general grounds to have Poincar\'e duality. 

Next, following the philosophy of the first paragraph, we are now led in Section~\ref{secAd} to introduce and study  the finite-dimensional quotient Hopf algebras
\[ A_d:=\F_p[x]/(g_d),\quad d\ge 1\]
which for regular $m$ of degree $d$ now inherit a connected calculus with $H^0_{\rm dR}(A_d;m)=\F_p1$. These algebras are much bigger than $\F_p[x]/(x^p)$ and moreover they fit into a directed system of Hopf algebras
\[ \{A_j\twoheadrightarrow A_i\ |\ i\ {\rm divides}\ j\}\] 
ordered by divisibility.  Here the polynomial $g_d$ is known from Artin-Schreier theory \cite[Ch. VI/Thm 6.4]{Lang} to be the product of all monic irreducibles of degree dividing $d$, so that $g_i$ divides  $g_j$ whenever $i$ divides $j$, leading to the map stated. The inverse limit $  \widehat{\F_p[x]}:=\lim_\leftarrow A_d$  
projects on to every $A_d$ and comes with a map $\F_p[x]\to \widehat{\F_p[x]}$ through which the quotienting maps $\F_p[x]\twoheadrightarrow A_d$ necessarily factor. Since each monic irreducible gives a field extension, it is also clear that  
\[ A_d\cong\prod_{k|d}\F_{p^k}^{N_k},\quad \widehat{\F_p[x]}=\prod_{m} {\F_p[x]\over (m) }\cong\prod_{k}\F_{p^k}^{N_k}\] as rings, where $N_k$ is the number of monic irreducibles of degree $k$. 
The coproducts imply that $\widehat{\F_p[x]}$ has some form of limiting Hopf algebra structure, but not with respect to the algebraic tensor product. This limit and its geometry are beyond our scope at present, where we  focus on the structure and geometry of the $A_d$ individually while thinking of them only loosely as increasingly good approximations of $\F_p[x]$.  Here one of the maps in the directed system is $A_d\twoheadrightarrow A_1$ for all $d$ and our main result (Theorem~\ref{cocycle}) is a structure theorem that  $A_d=C_d\tens_\chi A_1$ as a Hopf algebra cocycle extension, where $C_d$ is a sub-Hopf algebra generated by $g_1$ with a single relation related to the trace map.  This also implies that the $A_d$ are (cleft) Hopf-Galois extensions or (trivial) quantum principal bundles in the sense explained in \cite{Ma:ltcc,Ma:prim}. 

The paper concludes in Section~5 with a more detailed study of $A_1$ for general $p$ and $A_2$ for $p=2$. The geometric picture here is as the algebra of functions on $\Z/p\Z$ in the first case and a cocycle extension of the algebra of functions on $(\Z/2\Z)^2$ in the second. We focus on two important aspects; one is to compute the Hopf algebra Fourier transform and the other is to compute the moduli of translation-invariant `quantum metrics' and their associated quantum Levi-Civita connections. The latter turn out all to  be flat, which is consistent with the geometric picture but does take some proof in the case of $A_2$. Our results in this section are limited, but they suggest that $A_d$ in general could be an interesting class of finite-dimensional noncommutative geometries for further study.

\section{Preliminaries}\label{secpre} 

We will need the notion of a Hopf algebra $(A,\Delta,\eps,S)$ over a field $k$, where $A$ is a unital algebra and $\Delta a=a\o\tens a\t$ in a compact `Sweedler notation' (a sum of such terms understood) is an algebra homomorphism which, together with the counit character $\eps:A\to k$,  forms a coalgebra (or makes $A^*$ into an unital algebra). In addition, we require an antipode $S:A\to A$ obeying $(Sa\o)a\t=1\eps(a)=a\o Sa\t$ for all $a\in A$. More details can be found in many texts, including \cite{Ma:book}. We now give some preliminaries on noncommutative or `quantum' differentials central to the paper. 

\subsection{Noncommutative differentials} By definition, a {\em first order differential calculus} on a unital algebra $A$ is a pair $(\Omega^1,\extd)$ where $\Omega^1$ is an $A-A$-bimodule and $\extd:A\to\Omega^1$ obeys the Leibniz rule. We also require that $A\tens A\to \Omega^1$ given by sending $a\tens b$ to $a\extd b$ is surjective (otherwise one has a `generalised differential structure'\cite{MaTao}). The calculus is {\em connected} if $\ker\extd=k.1$ and {\em inner} if there exists $\theta\in\Omega^1$ such that $[\theta,a]=\extd a$ for all $a\in A$. As mentioned in the introduction, there is a universal first order differential calculus built on $\ker(\cdot)\subset A\tens A$ with $\extd a=1\tens a-a\tens 1$ which when $A$ is commutative quotients down to the K\"ahler differential. For higher forms we use the notion of a DGA over $A$,  meaning a graded algebra   $\Omega(A)=\oplus_n\Omega^n$ where $\Omega^0=A$, equipped with a graded-derivation $\extd:\Omega^i\to \Omega^{i+1}$ with respect to the product $\wedge$ of $\Omega$ and obeying $\extd^2=0$. More specifically, we require $\Omega$ to be generated by $\Omega^1$ and $A$ for a specified first order differential calculus, which is more restrictive than a regular DGA in other contexts (one says that $\Omega$ is the exterior algebra of the first order calculus). From this point of view, connectedness means $H_{\rm dR}^0=k.1$ as part of the cohomology of $(\Omega,\extd)$. For an inner DGA, we require $\extd=[\theta,\ \}$, a graded commutator. 

When $A$ is a Hopf algebra, a first order differential calculus is left covariant if $\Delta_L(a\extd b)=a\o b\o\tens a\t\extd b\t$ is well-defined as a map $\Omega^1\to A\tens\Omega^1$. If so, it becomes an $A$-coaction and $\Omega^1$ a left Hopf module. It follows from Hopf algebra theory that $\Omega^1$ is freely generated over $A$ by its space $\Lambda^1$ of left-invariant 1-forms. One can show that these are the image of the map $\varpi:A^+\to \Omega$ defined by $\varpi(a)=Sa\o\extd a\t$  (the `Maurer-Cartan form'), where $A^+=\ker(\eps:A\to k)$, with the result that $\Lambda^1\cong A^+/\CI$ for some right ideal $\CI$. Hence, classifying left covariant calculi amounts to classifying ideals in $A^+$. Moreover, the calculus is bicovariant (i.e., both left and right covariant) if and only if this ideal is also stable under a suitable adjoint coaction. In this case $\Lambda^1$ becomes a right $A$-crossed module (or Drinfeld-Yetter module) which in turn leads to a canonical construction for a braided exterior algebra $\Lambda$ via the braiding of this category when $S$ is invertible. Specifying the semidirect product of elements of $A$ with invariant forms then determines the full structure of the canonical extension $(\Omega,\extd)$ in the bicovariant case. The theory goes back to \cite{Wor} with a  treatment more along the above lines in \cite{Ma:prim,Ma:ltcc}. We omit details here since in our examples $A$ will be cocommutative, the adjoint coaction trivial (so all left-covariant calculi are bicovariant) and the braiding the trivial transposition map so that  $\Lambda$ is just the usual exterior algebra of $\Lambda^1$. In characteristic $2$, it means that all elements of $\Lambda^1$ square to zero.  

In view of this general picture, it is enough for a left-covariant or bicovariant calculus to focus on constructing and classifying the choice of $\Omega^1$. Here  a morphism between calculi over a fixed algebra $A$ means a bimodule map forming a commuting triangle with $\extd$ (this is a special case of the notion of a differentiable map between algebras equipped with first order differentials). Given the theory above, translation-invariant calculi on $k[x]$ are given by ideals in  $k[x]^+=(x)$ (the polynomials with no constant term) and irreducible differential calculi (those with no proper quotients) are given by $\Lambda^1=(x)/(x m(x))$ where $m$ is a monic irreducible polynomial \cite{Ma:fie,Ma:prim}. Clearly, $m(x)$ also defines a field extension $K=k[\mu]/(m(\mu))$ and as a vector space over $k$ one can identify $\Omega^1=K[x]$,  with $\mu^0=\extd x,\mu^1,\cdots,\mu^{d-1}$ a natural basis of $\Lambda^1$ over $k$ and of $\Omega^1$ over $k[x]$. Differentials do not commute with functions  except when $m(x)=x$, which is the classical calculus.  This leads to the explicit bimodule relations and $\extd$ stated in the  introduction, which look very much like a finite-difference calculus but should be interpreted as above with $\mu\in K$ a 1-form, not a parameter. For $d=1$, however, one can identify $\mu$ with an element of $k$, including 0 for the classical calculus. Other than the classical case, the calculus is inner with $\theta=\mu^{-1}$ computed in $K$.

\subsection{Construction of connected caluculi}

We will need a couple of observations relevant to the paper. The first applies to $A$ augmented by a morphism of unital algebras $\eps:A\to k$ (so $\eps(1)=1$) as a `base-point'. In the Hopf algebra case we will use the counit. Also in the Hopf algebra case, if $\Omega^1$ is bicovariant then its canonical  $\Omega$ is a super-Hopf algebra with $\Delta|_{\Omega^1}=\Delta_L+\Delta_R$ \cite{Brz} and we say more generally that $\Omega$ is strongly bicovariant when this happens. Here a super-Hopf algebra is like a Hopf algebra but we use the graded tensor product algebra, with respect to which $\Delta$ is a homomorphism. It was shown in \cite{MaTao} that $\extd$ is then also a super-coderivation. 

\begin{proposition} Let $(A,\eps)$ be an augmented unital algebra and $\Omega(A)$ an exterior algebra DGA. Then $A_c=A/\<H_{\rm dR}^0(A)\cap\ker\eps\>$ acquires an inherited differential calculus. If $A$ is a Hopf algebra and $\Omega(A)$  strongly bicovariant then $H_{\rm dR}^0(A)$ is a sub-Hopf algebra, $A_c$ is a quotient Hopf algebra and  the inherited $\Omega(A_c)$ is strongly bicovariant.  \end{proposition}
\proof (i) Clearly $J=H_{\rm dR}^0(A)\cap\ker\eps$ is a subalgebra by the Leibniz rule and we quotient by the ideal it generates intersected with the ideal $A^+$. More generally, we define $\Omega(A_c)=\Omega(A)/\<H_{\rm dR}^0(A)\cap\ker\eps\>$ where we quotient by the ideal generated in the exterior algebra. The map $\extd$ descends to this quotient by definition, for example $\extd (g a)=(\extd g)a+g\extd a=g\extd a$ for all $g\in J$. (ii) For the second part, by assumption there exist left and right coactions $\Delta_L,\Delta_R$ commuting with $\extd$. Hence if $a\in \ker \extd=H_{\rm dR}^0$ it follows that $(\id\tens\extd)\Delta a=\Delta_L\extd a=0$ and $(\extd\tens\id)\Delta a=\Delta_R\extd a=0$. Hence $\Delta a\in H_{\rm dR}^0\tens H_{\rm dR}^0$ and we have a sub-Hopf algebra. It follows that $J=H_{\rm dR}^0{}^+$ is a coideal (here if $a\in J$ then $\Delta a=a\o\tens (a\t-1\eps(a\t))+a\tens 1\in A\tens J+J\tens A$). In this case $I=\<J\>$, the ideal it generates, is a Hopf ideal and $A_c=A/I$ is a Hopf algebra. Employing the construction above, the assumed $\Delta_{L,R}$ descend to $A_c$ since these are part of Hopf module structures, e.g. $\Delta_L(\omega a)=(\Delta_L\omega)\Delta a\subseteq AJ\tens\Omega^1+A\tens\Omega^1 J$ for $a\in J$ and $\omega\in \Omega^1$. Hence  $\Omega^1(A_c)$ is bicovariant. The higher order calculi will be generated by $\Omega^1,A_c$ with the inherited relations and since the latter are bicovariant, the higher forms will be too. Equivalently, $\Omega(A)$ is a super-Hopf algebra and we can view $J$ as generating a super-coideal. \endproof

It is not clear that the quotient DGA now has $H_{\rm dR}^0(A_c)=k1$ but this is a step in the right direction and could be iterated. However, in our application  this construction will do the job in one step. We will also need the following lemma.

\begin{lemma}\label{H0g} Suppose that $H_{\rm dR}^0(k[x];m)\ne k.1$ and let $g\in H_{\rm dR}^0$ have minimal positive degree. Then $H_{\rm dR}^0(k[x];m)=k[g]$. \end{lemma}
\proof Let $f\in H_{\rm dR}^0$ and let $f=f_1 g+r_1$ where $r_1=0$ or $\deg r_1<\deg g$. Then $\extd r_1+(\extd f_1)g=0$ as $\extd g=0$. Viewing this
in the ring $K[x]$ we have the first term degree less than $ \deg g$ and the 2nd term degree greater than or equal to $\deg g$, hence both terms vanish and since $g$ had minimal degree among non-constants in $H_{\rm dR}^0$ we conclude that $r_1$ is a constant, $\extd f_1=0$. Iterating, we conclude that $f$ is a polynomial in $g$. \endproof

\section{ Structure of $H_{\rm dR}^0(\F_p[x]; m)$}

We consider $A=\F_p[x]$ and calculi defined by monic irreducible $m$ of degree $d$ and let $P$ denote the linear subspace of $\F_p[x]$ involving only power-$p$ exponents (we will explain later that $P$ is the set of primitive elements of $A$ as a Hopf algebra). If  $f\in P$ then $f$ is additive and hence, when $\mu\ne0$,
\[ \extd f=f(\mu)\mu^{-1}\in\F_{p^d}, \quad [v,f]=f(\mu)v,\quad\forall v\in \F_{p^d}\]
where $\F_{p^d}\subset\F_{p^d}[x]$ is the subspace of left-invariant 1-form in the calculus. Thus $f\in H_{\rm dR}^0\cap P$ are characterised by $f(\mu)=0$, while general $f\in H_{\rm dR}^0$ are characterised by $f(x+\mu)=f(x)$ which implies $f(\mu)=0$ for $f\in H_{\rm dR}^0\cap\ker\eps$ as a necessary condition. 

The case $d=1$ is easy enough to analyse in full detail and includes the case where $\mu=0$. Here monic irreducibles have the form $m(x)=x-\mu$, for some $\mu\in \F_p$. This gives a 1-dimensional calculus and of course no actual extension of the field. The field extension defined by $m$ would have generator set equal to $\mu$, which is why we have denoted this as the constant to fit with our previous notation. The case $\mu=0$ also leads to a calculus which we understand as the classical one. We let
\[ g_1(x)=x^p-x=x(x^{p-1}-1)=x(x-1)(x-2)\cdots(x-(p-1))\in P.\]

\begin{proposition}\label{H0d1} For $d=1$, $ H_{\rm dR}^0(\F_p[x];m)=\begin{cases} \F_p[x^p] &{\rm if}\  \mu=0 \\ \F_p[g_1] &{\rm if}\ \mu\ne 0\end{cases}$. \end{proposition}
\proof If $\mu=0$ we have the classical calculus where $\extd x^m=m x^{m-1}\extd x=0$ when $m=p$, and clearly any non-constant polynomial of lower degree will not be in the kernel by looking at its top degree. We then use Lemma~\ref{H0g}. When $\mu\ne 0$ we manifestly have $g_1(x+\mu)=g_1(x)$ (from the form of $g_1$) and we  show that its degree $p$ is the minimal degree of non-constant elements of $H_{\rm dR}^0$. Thus, let $f$ be monic of degree $t<p$ so $f=x^t+c x^{t-1}+\cdots$ for some $c\in \F_p$. We have $f(x+\mu)=x^t+ \mu t x^{t-1}+ c x^{t-1}+\cdots$ where we indicate further terms of degree less than $t-1$. For this to equal $f$ we need $\mu t=0$ mod $p$, which requires $t=0$ as $t<p$. Hence $f=1$. We then use Lemma~\ref{H0g}. \endproof

More generally, 
\begin{equation}\label{gn}  g_n(x):= x^{p^n}-x\in P\end{equation}
is, as mentioned in the introduction, the product of all irreducible monics in $\F_p[x]$ of degree dividing $n$. We are interested in a fixed monic $m$ of degree $d$ defining our differential calculus and associated $\F_{p^d}=\F_p[\mu]/(m)$. 

\begin{lemma}\label{gd} For $m$ of degree $d>1$, we have $H_{\rm dR}^0(\F_p[x]; m)\supseteq \F_p[g_d]$ and $g_i\notin H_{\rm dR}^0(\F_p[x];m)$ for $i=1,2,\cdots,d-1$. 
\end{lemma}
\proof Clearly $m$ is a factor of $g_d$ so $g_d(\mu)=0$. This is also immediate from $\mu^{p^d-1}=1$ in $\F_{p^d}$. Hence  we have $g_d\in H_{\rm dR}^0(\F_p[x];m)$ and hence (by the Leibniz rule) that all polynomials of it are contained in the cohomology. If $g_i(\mu)=0$ for some $i<d$ then $\mu$ is a zero of some irreducible monic of degree dividing $i$ and hence of degree less than $d$. This would have to be divisible by $m$, which is a contradiction. Hence $\extd g_i\ne 0$ for $1\le i<d$. \endproof

We say that $m$ of degree $d$ is {\em regular} if $H_{\rm dR}^0(\F_p[x];m)=\F_p[g_d]$. We have seen that this happens for $d=1$ precisely when $\mu\ne 0$. We will also be interested in 
\begin{equation}\label{hd} h_d(x):=x^{p^{d-1}}+x^{p^{d-2}}+\cdots+x\in P\end{equation}
where $h_d$ is the trace for the field extension when viewed as a map $\F_{p^d}=\F_p[\mu]/(m)\to \F_p$. 

\begin{lemma}\label{hd0} If $h_d(\mu)= 0$ then $m$ of degree $d$ is not regular, and $H_{\rm dR}^0(\F_p[x]; m)\supseteq \F_p[h_d]$ when $d>1$.
\end{lemma}
\proof If $d=1$ and $h_1(\mu)= 0$ then $\mu=0$ and we know that this case is not regular by the above. If $d>1$ and $h_d(\mu)=0$ then $h_d\in H_{\rm dR}^0$ by the above remarks and hence so is the subalgebra $\F_p[h_d]$ by the Leibniz rule. We also have
\begin{equation}\label{gdhd} g_d=h_d(h_d^{p-1}-1)\end{equation}
so that $\F_p[g_d]\subsetneq\F_p[h_d]$ and this is strict as $h_d$ has lower degree and clearly can't be written as a polynomial in $g_d$. Hence if $h_d(\mu)=0$ then $m$ cannot be regular. \endproof

\begin{theorem}\label{H0d} $m$ of degree $d$ has $H_{\rm dR}^0(\F_p[x];m)=\F_p[g_d]$, i.e. is regular, if and only if $h_d(\mu)\ne 0$.  
\end{theorem}
\proof $d=1$ was already covered so we fix $d\ge 2$ and prove the following assertion for $h_d(\mu)\ne 0$ by induction for $n$ in the range $1\le n\le d$: if $f$ is of degree less than $p^{n}$ and $f(x+\mu)-f(x)=c(\mu)$ for some $c$ in the $\F_p$-span of $\{\mu^{p^{i}}\ |\ n\le i\le d-1\}$  then $f(x)=f(0)$ and $c=0$. (The spanning set is empty when $d=n$ so $c=0$ in this case.) We note that if $f$ has degree less than $p$ and obeys the condition stated with $n=1$ then the argument in the proof of Proposition~\ref{H0d1} is unaffected when we look at powers of $x$   $>1$ and similarly allows us to conclude that $f$ has degree at most 1, so $f(x)=f(0)+a x$. Then  $f(x+\mu)-f(x)=c$ is $a\mu=c$ which requires $a=c=0$ as $\mu$ generates a normal basis when $h_d(\mu)\ne 0$. Thus the inductive assertion holds for $n=1$.  

Now suppose the assertion holds for $n-1$ in place of $n$ above and consider $f$ of degree less than $p^n$  and write this as $f=\sum_{k=0}^{p-1} x^{p^{n-1}k}f_k(x)$ where
$f_k$ have degree less than $p^{n-1}$. We also write
\begin{eqnarray*} f(x+\mu)&=& \sum_{k=0}^{p-1}(x^{p^{n-1}}+\mu^{p^{n-1}})^kf_k(x)+ \sum_{k=0}^{p-1}(x^{p^{n-1}}+\mu^{p^{n-1}})^{k}(f_k(x+\mu)-f_k(x))\\
&=&f(x)+A_{p-1}(x)\\ 
A_{r}(x)&=& \sum_{k=0}^{r}x^{p^{n-1}k}(f_k(x+\mu)-f_k(x))+\sum_{k=1}^{r}\sum_{s=0}^{k-1}\mu^{p^{n-1}(k-s)}x^{p^{n-1}s} \left({k\atop s}\right) f_k(x+\mu)\\
&=&x^{p^{n-1}r}(f_r(x+\mu)-f_r(x))+\sum_{s=0}^{r-1}\mu^{p^{n-1}(r-s)}x^{p^{n-1}s} \left({r\atop s}\right) f_r(x+\mu)+A_{r-1}.
\end{eqnarray*}
Now suppose that $f(x+\mu)=f(x)+c$ for $c$ in the span as asserted, i.e., $A_{p-1}(x)=c$. We prove that this implies that $c=0$ and $f$ is constant. Indeed, suppose $A_r(x)=c$. From the second expression for $A_r(x)$, only the first term has powers of degree greater than or equal to $p^{n-1}r$ which  tells us that $f_r(x+\mu)-f_r(x)=0$  and hence by our inductive hypothesis,  $f_r(x)=f_r(0)$ is a constant. Putting in this information gives us
\[ A_r(x)=\sum_{s=0}^{r-1}\mu^{p^{n-1}(r-s)}x^{p^{n-1}s} \left({r\atop s}\right) f_{r}(0)+A_{r-1}(x). \]
We now pick off the $\ge p^{n-1}(r-1)$ degrees to find for $r>1$,
\[ f_{r-1}(x+\mu)-f_{r-1}(x)+r\mu^{p^{n-1}} f_{r}(0)=0\]
and our induction hypothesis allows us to conclude that $f_{r}(0)=0$ and hence that $A_{r-1}(x)=c$ also. Starting at $r=p-1$ we now iterate this argument to conclude that
$f_{p-1}=0,\cdots,f_2=0$, $f_1$ is constant and $f_0(x+\mu)-f_0(x)+ \mu^{p^{n-1}}f_1(0)=c$, and hence, by our inductive hypothesis, that $f_0$ is constant and $f_1=c=0$. Hence our assertion also holds for $n$, completing the proof by induction. 

In particular, we apply this result with $n=d$ and $c=0$ to conclude that  $H_{\rm dR}^0$ contains no nonconstant elements of degree less than $p^d$. Hence the degree $p^d$ of $g_d$ is minimal among nonconstants in $H_{\rm dR}^0$. We then use Lemma~\ref{H0g}. Lemma~\ref{hd0} provides the other direction when $h_d(\mu)=0$. \endproof


\begin{corollary}\label{numreg}  $h_d(\mu)\ne 0$  iff  $m$ of degree $d$  has a nonzero coefficient in degree $d-1$.   Moreover, there are \[ {p-1\over pd}\sum_{k|d;  p\nmid k}\mu_\text{M\"ob}(k)p^{d\over k}\]
such $m$, where $\mu_\text{M\"ob}$ is the M\"obius function.
\end{corollary} 
\proof Here $h_d(\mu)=\Tr(\mu)$ is the trace for $\F_p[\mu]/(m)\to \F_p$ and it is a fact from number theory\cite[Ch. VI/Thm.~5.1]{Lang} that $\Tr(\mu)=-m_{d-1}$ where  $m(x)=x^d+m_{d-1}x^{d-1}+\cdots+m_0$ is the minimal polynomial of $\mu$, which is our case by construction. Hence $h_d(\mu)\ne 0$ if and only if $m_{d-1}\ne 0$. Next, the number of monic irreducibles in $\F_p[x]$ with a fixed non-zero value of this coefficient was found by Carlitz\cite{Car} and more recently in the form we use in \cite{Rus}. As we required only a non-zero value, we multiply this by the $p-1$ possible values to give the expression stated. \endproof

This gives an easy criterion to tell if a given $m$ is regular.  The number of such should be compared with Gauss' formula for the number $N_d$ of all irreducible $m$ of degree $d$,  $N_d={1\over d}\sum_{k|d}\mu_\text{M\"ob}(k)p^{d\over k}$. Thus a good fraction of $m$ are regular.  The formula gives $p-1$ regular $m$ as it should for $d=1$ by Proposition~\ref{H0d1}. 

Also note that the factorisation (\ref{gdhd}) means that either $h_d(\mu)=0$ (the non-regular case) or $m$ divides $(h_d-1)(1+h_d+\cdots h_d^{p-2})$ (the regular case) according to our theorem. Meanwhile, Lemma~\ref{hd0}  suggests a similar result for the cohomology for the flip side when $h_d(\mu)=0$:

\begin{conjecture}\label{H0conj}  $m$ of degree $d>1$ has  $H_{\rm dR}^0(\F_p[x];m)=\F_p[h_d]$ if and only if $h_d(\mu)=0$. 
\end{conjecture}

It is not clear that this can be proven by similar methods to those of our main theorem. We also note in passing that as well as the Trace there is a norm map $N:\F_{p^d}\to \F_p$ defined as 
\[ N(x)=x x^p\cdots x^{p^{d-1}}=x^{1+p+p^2+\cdots p^{d-1}}=x^{[d]_p},\quad [d]_p={p^d-1\over p-1}\]
and  $N(\mu)=(-1)^dm_0\ne 0$ as $m$ is irreducible. Hence $N\notin H_{\rm dR}^0(\F_p[x];m)$. 

\begin{example}\label{Exp2} Conjecture~\ref{H0conj} is supported by computer calculations   for $p=2$ and $d\le 4$ with code available on \cite{GIT} and with the following verified up to polynomials of degree 100:
\begin{enumerate}\item $H_{\rm dR}^0(\F_2[x];\mu^2+\mu+1)=\F_2[x^4+x]=\F_2[g_2]$
\item $H_{\rm dR}^0(\F_2[x];\mu^3+\mu^2+1)=\F_2[x^8+x]=\F_2[g_3]$
\item $H_{\rm dR}^0(\F_2[x];\mu^3+\mu+1)=\F_2[x^4+x^2+x]=\F_2[h_3]$
\item $H_{\rm dR}^0(\F_2[x];\mu^4+\mu^3+\mu^2+\mu+1)=\F_2[x^{16}+x]=\F_2[g_4]$
\item $H_{\rm dR}^0(\F_2[x];\mu^4+\mu^3+1)=\F_2[x^{16}+x]=\F_2[g_4]$
\item $H_{\rm dR}^0(\F_2[x];\mu^4+\mu+1)=\F_2[x^8+x^4+x^2+x]=\F_2[h_4]$
\end{enumerate}
where $h_3(\mu)=\mu(\mu^3+\mu+1)=0$ in (3) and $h_4(\mu)=\mu(\mu^3+1)(\mu^4+\mu+1)=0$ in (6) and $h_d(\mu)\ne 0$ in the other cases. This also illustrates Theorem~\ref{H0d} and  Corollary~\ref{numreg}, now proven. One can moreover see here that the regular $m$ are precisely the factors of degree $d$ in $h_d+1$, as per the general theory when $p=2$. 
\end{example}

\section{The Hopf algebras $A_d$}\label{secAd}

Motivated by the above cohomology computations, for each $d\in \N$ and each regular $m$ of degree $d$, we define
\[ A_d:=\F_p[x]/( g_d)=\F_p[x]/ ( x^{p^d}-x),\quad  \Omega(A_d;m):=\Omega(\F_p[x]; m)/\<g_d\>.\]
Here $J= H_{\rm dR}^0\cap \ker\eps={\rm span}\{g_d^m\ |\ m>0\}=\F_p[g_d]^+$ where the $+$ denotes functions with no constant term. Hence $\F_p[x]J=(g_d)$ is the ideal that we quotient out by to define $A_d$.

We think of the algebra $A_d$ as defining a `space' at a topological level in some sense, and in this regard we note that $A_d$ as defined depends only on the degree $d$ of the field extension. We think of $m$ as adding to this data a differentiable structure inherited from the one on $\F_p[x]$ or equivalently an element $\mu$ of a field extension of degree $d$.

\begin{corollary}\label{H0Ad} $H_{\rm dR}^0(A_d;m)=\F_p1$ for the inherited differential structure from any regular monic irreducible $m$ of degree $d$.
\end{corollary}
\proof This is immediate from Theorem~\ref{H0d} and we use the notations there. Suppose that $f(x+\mu)-f(x)\in (g_d)$ in $\F_{p^d}[x]$, i.e. products of $g_d$ with polynomials that include powers of $\mu$ in their coefficients.  In $\F_p[x]$ we let $f=h g_d+r$ where  either $r=0$ (so $f=0$ in $A_d$) or the degree of $r$ is less than $p^d$. Then $h(x+\mu)g_d(x+\mu)+r(x+\mu)-h(x)g_d(x)-r(x)=(h(x+\mu)-h(x))g_d(x)+r(x+\mu)-r(x)\in (g_d)$ since $g_d(x+\mu)=g_d(x)$ as $g_d$ is additive and $g_d(\mu)=0$. Hence $r(x+\mu)-r(x)\in (g_d)$. But since every nonzero element of $(g_d)$ has degree greater than or equal to $p^d$ we conclude that $r(x+\mu)-r(x)=0$ and hence by Theorem~\ref{H0d} that $r$ is a constant, hence $f$ is a multiple of the identity in $A_d$. .  \endproof

This means that we achieved our goal of having finite-dimensional quotients of $\F_p[x]$ equipped now with (a moduli space of) connected differential calculi. We now turn to the algebraic structure of the $A_d$. 
 
 \begin{corollary} $A_d$ is a $p^d$-dimensional Hopf algebra and $\Omega(A_d)$ is bicovariant. Moreover,  the primitive elements of $A_d$ are spanned by the set $\{x^{p^i}  : 0 \leq i < d \}$
  \end{corollary}
\proof As $x$ is primitive, we have
$$\Delta x^n = \sum_{k=0}^{n} \binom{n}{k} x^k \otimes x^{n-k}.$$
By Lucas' theorem\cite{Lucas,Fine}, $\binom{n}{k} = 0$ mod $p$ iff a base $p$ digit of $k$ is greater than the corresponding digit of $n$.  Hence $x^n$ is primitive if $n = p^i$ for some $i$.  Conversely, if $n$ is not a power of $p$, then there exist some $k\ne 0,n$ for which $\binom{n}{k} \neq 0$ mod $p$, and so $x^n$ is not primitive.  It then follows easily that the primitive elements of $\F_p[x]$ are precisely spanned by the $p$-power exponents.  In particular, $g_d$ is primitive and hence $A_d$ is a Hopf algebra. Its primitives have the same form but restricted to degree less than $p^d$. Here $\F_p[g_d]$ is a Hopf algebra with $g_d$ primitive and in this way a sub-Hopf algebra of $\F_p[x]$ as per the general theory in Section~2. The latter   also implies bicovariance. \endproof

Next, $A_d$ carries the Frobenius automorphism $F(x)=x^p$ and hence always contains a Frobenius-fixed subalgebra $B_d=A_d^F$ where every element equals its $p$-power. For $p=2$ this means that $B_d$ is a Boolean subalgebra.  
\begin{proposition}\label{dimBd}
$ \text{dim}\; B_d = \frac{1}{d} \sum_{k | d} \phi(k) p^\frac{d}{k}$, the number of irreducible factors of $g_d$. Here $\phi$ is the Euler totient function.
\end{proposition}
\proof  The Frobenius automorphism has order $d$ and permutes the set $\{1, x, x^2, \ldots, x^{p^d-1}\}$.  Write this permutation in its decomposition as cycles $\sigma_1, \ldots \sigma_b$.  When a polynomial $f$ is the sum of monomials from an orbit of a $\sigma_i$ it is fixed by the endomorphism.  The set of such polynomials (with all coefficients $1$) is linearly independent and generates $B_d$. Now let $C_s = \{sp^j \mod p^d - 1 : 0 \leq j \leq d-1\}$ be the cyclotomic coset of $p$ modulo $p^d - 1$ containing $s$.  
Note that each $C_s$ and $C_r$ are either disjoint or equal.  
Let $\mathcal{C} \subset \Z/(p^d -1)\Z$ be such that 
$$\bigcup_{s\in \mathcal{C}} C_s = \Z/(p^d-1)\Z$$ 
and each pair $C_s$ and $C_r$ are disjoint for $s,r\in\mathcal{C}, \; s \neq r$. 
$\mathcal{C}$ is in bijection with the set of orbits of the permutation cycles define above,
excluded the singleton orbit $\{x^{p^d -1} \}$, which comes from the additional factor of $x$ in
the polynomial modulus: let $s \in\mathcal{C}$, if $x^s \in \text{orb}\sigma_i$, then 
$\text{orb} \sigma_i = \{ x^{sp^j} \mod x^{p^d} - x: 0 \leq j \leq d-1 \}$.  

Let $\alpha \in \mathbb{F}_{p^d}$ be a generator for the multiplicative group $\mathbb{F}_{p^d}^\times$.  It is known \cite[Thm.~3.4.11]{LiXi} that $$x^{p^d-1} -1 = \prod_{s \in \mathcal{C}} m_s(x) $$ for $m_s(x) = \prod_{a\in C_s} (x-\alpha^a)$, hence the set of orbits of the permutation, and thus the basis we have given, is in bijective correspondence with irreducible factors of $x^{p^d}-x$.\endproof

Next recall from the introduction that as part of the inductive system there are canonical Hopf algebra maps $\pi:A_d\to A_1$, since $1$ divides every $d$.  To find the kernel of this map we consider the canonical automorphism of the algebra $A_d$ given by the order $p$ periodicity map $\psi(x)=x+1$. Here $g_i(x+1)=(x+1)^{p^i}-(x+1)=g_i(x)$ working over $\F_p$, so these give invariant elements of $A_d$ for $i=1,\cdots,d-1$. We will be interested in the invariant subalgebra $C_d=A_d^\psi$ of $A_d$.

\begin{proposition}\label{Cdext} $C_d=\F_p[g_1]/(h_d(g_1))$ for all $d\in\N$ is a Hopf algebra of dimension $p^{d-1}$ and 
\[ C_d\hookrightarrow A_d \twoheadrightarrow A_1\]
is an extension of Hopf algebras. 
\end{proposition}
\proof We start with $\F_p[x]^\psi=\F_p[g_1]$. This is known from Artin-Schreier theory but for completeness we include an elementary proof from \cite{Leon}. If $f(x+1)=f(x)$ and $f$ has degree less than $p$ then the proof of Proposition~\ref{H0d1} with $\mu=1$ applies and allows us to conclude that $f(x)$ is a constant. More generally let $f(x)=g_1 h(x)+r(x)$ where $r$ has degree less than $p$. Then $g_1(h(x+1)-h(x))=-(r(x+1)-r(x))$ which by degrees requires both $h$ and $r$ to be invariant. Thus $r$ is a constant and $h$ is an invariant of lower degree, leading to the result. Also clearly, the $\{g_1^i\}$ are linearly independent over $\F_p$ (by looking at the top degree of a polynomial relation). Also in $\F_p[x]$ we have $h_i(g_1)=g_1+g_1^p+\cdots g_1^{p^{i-1}}=x^p-x+(x^p-x)^p+\cdots+(x^p-x)^{p^{i-1}}=x^p-x+x^{p^2}-x^p+\cdots x^{p^i}-x^{p^{i-1}}=g_i$ on cancellation. In particular, $h_d(g_1)=g_d$, and if a polynomial in $g_1$ of degree less than $p^{d-1}$ is divisible by $g_d$ then, by degrees, it must separately vanish. Hence polynomials in $g_1$ up to degree less than $p^{d-1}$ viewed in $A_d$ form a $p^{d-1}$-dimensional subalgebra. Finally, if $f\in\F_p[x]$ has degree less than $p^d$ and $f(x+1)-f(x)$ is divisible by $g_d$ then by degrees is must separately vanish, hence $C_d=\F_p[g_1]/(g_d)=\F_p[g_1]/(h_d(g_1))$. This inclusion $i:C_d\hookrightarrow A_d$ makes $C_d$ a sub-Hopf algebra as $g_1$ is primitive. The Hopf algebra map $\pi:A_d\to A_1$ where we quotient by $(g_1)$ clearly obeys $\pi \circ i=1\eps$ since it is $1$ on $1\in C_d$ and vanishes on $C_d^+$. It follows from general arguments since the Hopf algebras involved are finite dimensional, see \cite[Cor~3.2.2]{Sch}, that this gives an exact sequence of Hopf algebras in the technical (cleft) sense provided only that the dimensions match. This is our case as we have seen that $\dim(A_1)\dim(C_d)=\dim(A_d)$. \endproof

It follows from the theory of such extensions of Hopf algebras that $A_d$ is a cocycle bicrossproduct of $C_d$ and $A_1$ in the sense of \cite[Sec.~6.3]{Ma:book}. We will give this explicitly and find in fact that the extension result applies to $\F_p[x]$ as well, not only the finite-dimensional quotients.  We first define
\[ \delta_i(x)=-{g_1\over x-i}=-\prod_{j\ne i}(x-j)\in \F_p[x],\quad i=0,\cdots,p-1\]
which clearly obey $\psi(\delta_i)=\delta_{i-1}$. Hence $\sum_i\delta_i$ is $\psi$-invariant and has degree $p-1$ hence is a constant. Evaluating at zero, only $\delta_0(0)=-(p-1)!=1$ is non-zero, we have
\[  \sum_i\delta_i=1.\]
The $\delta_i(x)$ are similar to the $({x\atop i})$ basis functions in Mahler's theorem\cite{Mah} albeit the context is different. The following is presumably known but we have not found it elsewhere and include a short proof. 

\begin{lemma}\label{A1isom}  $A_1\cong \F_p(\Z/p\Z)$ the Hopf algebra of functions on the finite group $\Z/p\Z$. 
\end{lemma}
\proof We identify the Kronecker delta-function at $i\in \Z/p\Z$ with $\delta_i\in A_1$ i.e. viewed mod $g_1$.  Clearly, in $A_1$ we have $\delta_i(x)(x-i)=0$ and hence $\delta_i\delta_j=0$ in $A_1$ for $i\ne j$, so that $\delta_i\delta_i=\delta_i$ in $A_1$ from $\sum\delta_i=1$. Hence this is an isomorphism of algebras. For the coproduct we note that the image of the coproduct of $\F_p[x]$ has the property of invariance under $\psi\tens\psi^{-1}$ acting in the two factors (this clear for $\Delta x$ and therefore applies on any polynomial). Since $\{\delta_i\}$ by the above relations form a basis of $A_1$, we let $\Delta\delta_k=\sum_{i,j}c^k_{ij}\delta_i\tens\delta_j$ for some $c^k_{ij}\in \F_p$. Then invariance implies that $c^k_{i,j}=c^k_{0,i+j}$. However, $\eps(\delta_i)=\delta_{i,0}$ and the counity axiom then implies that $c^k_{0,i}=\delta_{k,i}$. Hence $\Delta\delta_i=\sum_{j=0}^{p-1}\delta_{i-j}\tens\delta_j$ as for $\F_p(\Z/p\Z)$. \endproof

\begin{theorem}\label{cocycle}  $\F_p[x]=\F_p[g_1]\tens_\chi A_1$  is a cocycle cleft extension by $A_1$ coacting via $\psi$ and for $p>2$ with cocycle 
\[ \chi:A_1\tens A_1\to \F_p[g_1],\quad \chi(\delta_i\tens\delta_j)=\begin{cases} 1 &\text{\rm if\ } i=j=0\\ {g_1\over j} &\text{\rm if\ }i=0,\ j\ne 0 \\ {g_1\over i} &\text{\rm if\ }i\ne0,\ j=0\\ -{g_1\over i}&\text{\rm if\ }i=j\ne 0\\ 0 &\text{\rm else.}\end{cases}\]
This amounts to the new identities for $\delta$-functions in $\F_p[x]$ for $p>2$,
\[ \delta_i\delta_i=\delta_i+g_1\sum_{k=1}^{p-1}{\delta_{i+k}\over k},\quad\delta_i\delta_j=-g_1{\delta_i-\delta_j\over i-j}\]
for all $i,j\in \F_p$ and $i\ne j$. For $p=2$ the cocycle and identities are
\[ \chi(\delta_i\tens\delta _j)=g_1+\delta_{i,0}\delta_{j,0},\quad \delta_i^2=\delta_i+g_1,\quad\delta_0\delta_1=g_1.\]
The coproduct of $\F_p[x]$ becomes
\[ \Delta \delta_i=\sum_{j=0}^{p-1}\delta_{i-j}\tens\delta_j,\quad \Delta g_1=g_1\tens 1+1\tens g_1,\quad \eps\delta_i=\delta_{i,0},\quad \eps g_1=0\]
so that $A_1$ is a subcoalgebra. These formulae descend to $A_d=C_d\tens_\chi A_1$ for all $d\ge 1$. 
\end{theorem}
\proof In view of Lemma~\ref{A1isom}, the action of $\Z/p\Z$ via $\psi$ on $\F_p[x]$ becomes a right coaction $\Delta_R f=\sum {\psi^i}(f)\tens\delta_i$ of $A_1$. 
Clearly $\Delta_R(\delta_i)=\sum_{j}\delta_{i-j}\tens\delta_j$ viewed in $\F_p[x]\tens A_1$. Then $\phi:A_1\to \F_p[x]$  sending $\phi(\delta_i)=\delta_i$ is a right comodule map. It is also convolution-invertible with $\phi^{-1}(\delta_j)=\delta_{-j}$ as 
\[ \sum_j \delta_j\delta_{j-i}=\delta_{i,0}.\]
This is because the sum is $\psi$-invariant hence by degrees is at most linear in $g_1$. The constant value is $\delta_{i,0}$ since only $\delta_0(0)=1$ is non-zero, while 
\[ \delta_0^2=1+O(x^2),\quad \delta_0\delta_j=-{1\over j}x+O(x^2),\quad \delta_i\delta_j=O(x^2)\]
for all $i,j\ne 0$. Using that $\sum_{i=1}^{p-1}1/i=0$ which is equivalent to $\sum_{i\in\F_p}i=0$ mod $p$ valid for $p>2$, one has that $\sum_j\delta_j\delta_{j-i}$ has zero coefficient in degree $1$, so there is no $g_1$ term. Hence we have a cleft extension and $\F_p[x]\cong \F_p[g_1]\tens_\chi A_1$ for some cocycle  $\chi:A_1\tens A_1\to \F_p[g_1]$ which we compute from
\[ \chi(\delta_i\tens\delta_j)=\sum_k \phi(\delta_{i-k})\phi(\delta_{j-k})\phi^{-1}(\delta_k)=\sum_k\delta_{i+k}\delta_{j+k}\delta_{k}.\]
This is again $\psi$-invariant and has degree at most $(p-1)^3$, so is at most quadratic in $g_1$. Looking to degree $2$, we have
\[ \delta_i\delta_j\delta_k=O(x^3),\quad \delta_i\delta_j\delta_0={x^2\over ij}+O(x^3),\quad \delta_i\delta_0^2=-{x\over i}-{x^2\over i^2}+O(x^3), \quad \delta_0^3=1+O(x^3)\]
where we used that $\sum_{i=1}^{p-1} 1/i^2=0$ mod $p$ for $p>3$, which is equivalent to  a power-sum identity $\sum_{i\in\F_p}i^2=0$ mod $p$ for $p>3$. Such identities are known to hold for all powers not divisible by $p-1$, see \cite{Carps}. From this one can see that $\chi$ has no $x^2$ term, and hence 
is at most a constant plus linear term in $g_1$. We then use our expression for $\chi$ and the form of triple products of $\delta$'s to match the constant terms and coefficients of $x$, giving the cocycle as stated. From the theory of extensions, see \cite[Prop.~6.3.2]{Ma:book}, we will  be able to recover the product of $\F_p[x]$ from
\begin{equation}\label{cocyprod}(c\tens \delta_i)(c'\tens \delta_j)=cc'\sum_{k}\chi(\delta_{i-k}\tens\delta_{j-k})\tens\delta_k,\quad\forall c,c'\in \F_p[g_1]\end{equation}
which implies in particular that 
\[ \delta_i\delta_j=\sum_{k=0}^{p-1}\chi(\delta_{i-k}\tens\delta_{j-k})\delta_k\]
holds in $\F_p[x]$. This provides the identities stated for $p>2$. For $p=2$ it is easy to verify the stated identities with $\delta_0=1+x$ and $\delta_1=x$ and $g_1=x^2+x$ in this case, from which the cocycle has to have the form stated. 

Finally, we look at the coproduct. Its image in $\F_p[x]\tens\F_p[x]$ is invariant under $\psi\tens\psi^{-1}$ which together with the factorisation as algebras already proven implies a general form $\Delta\delta_i=\sum_{j,k}\delta_jc^i_{j,k}\delta_k$ for $c^i_{j,k}\in\F_p[g_1]\tens\F_p[g_1]$. The same arguments as in the proof of Lemma~\ref{A1isom} apply and tell us that $\Delta\delta_i=\sum_{j,k}\delta_jc^i_{0,j+k}\delta_k$. Writing $c^i_{0,j}=\delta_{i,j}+g_1\tens 1 b^i{}_j+ 1\tens g_1 b'{}^i{}_j+(g_1\tens g_1)c'{}^i{}_j$ where $b,b'$ are constants, the counit axiom  $(\id\tens\eps)\Delta\delta_i=(\eps\tens\id)\Delta\delta_i=\delta_i$ tells us that $b=b'=0$ (and also fixed the first term as $\delta_{i,j}$). Now, $\Delta$ does not change the total degree so the total degree of 
\[ \Delta\delta_i=\sum_{j}\delta_{i-j}\tens\delta_j+ (g_1\tens g_1)\sum_{j,k}\delta_jc'{}^i_{j+k}\delta_k\]
has to be $p-1$. The first term has total degree at most $<2p$ while the second term has leading term $(x^p\tens x^p)\sum_{j,k}\delta_jc'{}^i_{j+k}\delta_k$, hence this second term must separately vanish.  This means that we have a tensor product as coalgebras, so that $A_1$ appears as a subcoalgebra. 

Clearly, these results are not changed modulo $g_d$ as higher powers of $g_1$ were not involved. Hence $A_d= C_d\tens_\chi A_1$ also by identifying the $\delta$-functions. \endproof

Such cleft extensions may also be regarded as trivial quantum principal bundles or Hopf-Galois extensions of a certain trivial type, which are indeed classified by cocycles as explained in detail in \cite[Sec.~5]{Ma:nonab}. Indeed, the above implies that they are of the quantum homogeneous space type with right coaction  $(\id\tens\pi)\Delta$ of the fibre Hopf algebra $A_1$ on the total space Hopf algebra $A_d$, with base algebra $C_d$. Geometrically by  Lemma~\ref{A1isom}, the underlying structure group is $\Z/p\Z$. 

\section{Noncommutative geometry of $A_1$ and $A_2$}

In this section we study $A_1,A_2$ in more detail, focussing on their Fourier theory and translation-invariant noncommutative differential geometry respectively. The aim is to obtain a fuller picture of these algebras as examples of the $A_d$ family. 

 Fourier transform works on any finite-dimensional Hopf algebra $A$ equipped with (say) a right translation-invariant integral $\int:A\to k$ and another $\int:A^*\to k$ on its dual. Here $(\int\tens\int)(\exp)\ne 0$, where $\exp=e_a\tens f^a$ is the canonical coevaluation elements for the duality pairing defined by any basis $\{e_a\}$ of $A$ with dual basis $\{f^a\}$.  Translation invariance means $(\int\tens\id)\Delta=1\int$. Fourier transform is then $\CF:A\to A^*$ defined by $\CF(f)=(\int e_a f)f^a$. The inverse is similar but with the Hopf algebra antipode, see \cite[Prop.~1.77]{Ma:book} for an exposition. Fourier transform is compatible with any translation-invariant differential calculus $\Omega^1$ and turns the translation-invariant differentials $\del^a$ with respect to a basis into right multiplication by the corresponding dual basis element $f^a$ in $A^*$. 

The formalism of noncommutative Riemannian geometry works for any algebra $A$ with differential structure defined at least to $\Omega^2$. By a `metric' we mean an element $g\in\Omega^1\tens_A\Omega^1$ which is quantum symmetric  in the sense $\wedge(g)=0$ and invertible in the sense of existence of a bimodule map $(\ ,\ ):\Omega^1\tens_A\Omega^1\to A$ such that $(\omega, g\uo)g\ut=\omega=g\uo(g\ut,\omega)$ for all $\omega\in\Omega^1$. Here $g=g\uo\tens g\ut$ (a sum of such terms understood) is a notation. One can show\cite{BegMa4} that such a $g$ is necessarily central.  By a `left connection', in our case on $\Omega^1$,  we mean $\nabla:\Omega^1\to \Omega^1\tens_A\Omega^1$ such that $\nabla(a\omega)=a\nabla\omega+\extd a\tens\omega$ for all $a\in A$ and $\omega\in \Omega^1$.  By a `bimodule connection'\cite{DV1,DV2} we mean a left connection such that in addition $\nabla(\omega a)=(\nabla\omega)a+\sigma(\omega\tens\extd a)$ for some bimodule map $\sigma:\Omega^1\tens_A\Omega^1\to \Omega^1\tens_A\Omega^1$. If a left connection admits such a $\sigma$ then the latter is unique, hence this is a property of $\nabla$ and not further data. In this case one has the notion of metric compatible connection $\nabla g=0$ where $\nabla$ acts on each tensor factor $\Omega^1$ and $\sigma$ is used to correctly position its output when acting on the second tensor factor. Finally, the torsion of a connection on $\Omega^1$ is  $T=\wedge\nabla-\extd:\Omega^1\to \Omega^2$ and in noncommutative Riemannian geometry we are ideally interested in finding a `Levi-Civita' bimodule connection defined as metric compatible and torsion free.  
More details can be found in \cite{Ma:ltcc}.

 We will be interested in the translation-invariant geometry in the case of $A$ a Hopf algebra. In practice this just means that the coefficients are constant with respect to basis of left-invariant 1-forms. Unlike Lie theory, the choice of invariant differential structure is not unique; we take the calculus on $A_d$ inherited from that of $\F_p[x]$ for a choice of regular $m$ of degree $d$.

\subsection{Fourier transform and geometry on $A_1$}

Here we focus on 
\[ A_1= \F_p[x]/(x^p-x)\]
as a Hopf algebra with $x$ primitive. This is isomorphic to $\F_p^p$ as a  ring, hence to the functions on $p$ points, and indeed we have already remarked in Lemma~\ref{A1isom} that it is isomorphic to $\F_p(\Z/p\Z)$ where the Kronecker delta-functions on the latter are mapped to the $\delta_i(x)$ for $i=0,\cdots,p-1$.  From the projector relations among the $\delta_i$ in $A_1$ and the evaluations $\delta_i(j)=\delta_{i,j}$
which follow from the definition of $\delta_i(x)$,  it is easy to see that 
\begin{equation}\label{deltaf} \delta_i(x)f(x)=f(i)\delta_i(x),\quad \forall f(x)\in A_1.\end{equation}
Thus, if $f(x)\in A_1$ then the values $f(i)$ for $i\in\Z/p\Z$ provide the corresponding function on the group while conversely  $f(x)=\sum_i\delta_i(x)f(i)$. We also recall that finite-dimensional Hopf algebras have unique translation-invariant integration up to normalisation. In our case we have up to normalisation
\begin{equation}\label{xint}  \int x^i=\begin{cases}1 & {\rm if}\ i=p-1\cr 0 &{\rm otherwise,}\end{cases}\end{equation}
which is equivalent via the isomorphism to  $\int f=\sum_{i=0}^{p-1} f(i)$ for $f\in \F_p(\Z/p\Z)$. From this or from the coefficient of $x^{p-1}$  in $\delta_i$ being 1, we clearly have $\int\delta_i(x)=1$. 

Next, the dual Hopf algebra to $\F_p(\Z/p\Z)$ is the group Hopf algebra of $\Z/p\Z$,
\[ {A_1^*}=\F_p \Z/p\Z= \F_p[t]/(t^p-1),\quad \Delta t=t\tens t\]
and has the unique normalised translation-invariant integral
\begin{equation}\label{tint} \int t^i= \begin{cases}1 & {\rm if}\ i=0\cr 0 &{\rm otherwise.}\end{cases}\end{equation}
We can view this Hopf algebra as dual to $A_1$ via the Hopf algebra duality pairing  
\begin{equation}\label{txpair}A_1\tens {A_1^*}\to \F_p,\quad \<f(x),t^j\>=f(j),\quad \forall f(x)\in A_1.\end{equation}
As our Hopf algebras are finite-dimensional, there is necessarily a canonical coevaluation which we denote  $\exp\in {A_1^*}\tens A_1$. We recall that for any finite dimensional Hopf algebra and specified integral on it, we have a Fourier transform $\CF:A_1\to {A_1^*}$ given by integration against one factor of $\exp$, see \cite[Prop.~1.7.7]{Ma:book} for an exposition. In our case it is immediate from (\ref{txpair}) that $\exp=\sum_{i=0}^{p-1} t^i \tens\delta_i(x)$ 
leading to the canonical Hopf algebra Fourier transform 
\begin{equation}\label{xtF}  \CF(f)=\sum_{i=0}^{p-1} t^i f(i),\quad \CF^{-1}(t^i)=\delta_i(x),\quad \forall f(x)\in A_1,\ i\in0,\cdots,p-1.\end{equation}
Note that $(\int\tens\int)\exp=1$. This completes our review of Fourier theory on $\Z/p\Z$. 

Next, in the same way as we have described the functions on the finite group as a quotient of the affine line $\F_p[x]$, namely $A_1$,  we can do the adjoint thing on the dual side. Thus, $\F_p\Z/p\Z$ already looks like an algebraic group with group-like generator $t$ but we can go further and write this as like the enveloping algebra of a Lie algebra with infinitesimal generator $L$, say.

\begin{lemma}\label{H1} Let $p>2$. 
\[ {A_1^*}= \F_p[L]/( L^p )\]
as a Hopf algebra via the identification
\[ t=e^L:=\sum_{i=0}^{p-1} {L^i\over i!},\quad  L=\ln(t):=-\sum_{i=1}^{p-1}{t^i\over i}\in {A_1^*}^+,\]
in terms of a `truncated exponential' $e^{(\ )}$ and `truncated logarithm' $\ln(\ )$. We have
\[ \int L^i=\begin{cases}1 &{\rm if\ }i=0, p-1\\ 0 &{\rm else}\end{cases} \]
as equivalent to (\ref{tint}). 
\end{lemma}
\proof First we note that given $t^p=1$ and $L$ defined as stated, $L^p=-\sum_{i\ne 0}1/i=0$ and 
\[ i L=i \ln(t)=\ln(t^i)   \]
for all integers $i$ mod $p$. Conversely, given $L$ with $L^p=0$, we define $t=e^L$ and clearly $t^p=\sum_{i= 0}^{p-1}(L^i)^p/i!=1$. More generally it follows from $L^p=0$  that
\[ e^{iL}e^{jL}=\sum_{k=0}^{p-1}\sum_{s=0}^{p-1}{i^k j^s\over k! s!}L^{k+s}=\sum_{m=0}^{p-1}\sum_{k=0}^m\left({m\atop k}\right){L^m\over m!} i^k j^{m-k}=e^{(i+j)L}\]
which implies in particular that $t^i=e^{iL}$ and hence 
\[ \ln(e^L)=-\sum_{i=1}^{p-1}{e^{iL}\over i}=-\sum_{i=1}^{p-1}{1\over i}\sum_{j=0}^{p-1}i^j {L^j\over j!}=-\sum_{j=0}^{p-1}{L^j\over j!}\sum_{i=1}^{p-1}i^{j-1}=-L(p-1)=L\]
using the power-sum identity so that the sum over $i$ contributes only for $j=1$. Hence the algebra map $\F_p[L]/(L^p)\to {A_1^*}$ sending $L$ to $\ln(t)$ is injective (by applying the algebra map going the other way that sends $t$ to $e^L$) and hence by dimensions an isomorphism. Next, given $L$, for $t$ to be group-like we need
\[ \Delta L=\ln(e^L\tens e^L),\quad\eps L=0.\]
With this coalgebra,  the two Hopf algebras are isomorphic. We then convert over the integral as stated. \endproof
 
We remark that the coproduct can be written more explicitly as
\begin{equation}\label{DeltaL}\Delta L=L\tens 1+1\tens L-\sum_{i=1}^{p-1} {\left({p\atop i}\right)\over p}L^i\tens L^{p-i}\end{equation}
which makes sense as $p$ divides the binomial coefficient. We have verified this by computer for small primes. The coproduct can also be written as a multiplicative correction
\begin{equation}\label{DeltaLnew}\Delta L=\left(1-\sum_{i=1}^{p-2}a_i L^i\tens L^{p-1-i}\right)(L\tens 1+1\tens L)\end{equation}
where 
\[  a_i={\left({p-1\atop i}\right)-(-1)^i\over p} \]
also make sense and obey $a_i=a_{p-1-i}$, $a_1=1$, $a_2=(p-3)/2$ etc, with middle value $i=(p-1)/2$ giving the so-called `swinging Wilson quotients'\cite{sw}.   Also note that if we took $L$ primitive then we would again have a Hopf algebra on $\F_p[L]/(L^p)$ but now dually paired with $\F_p[x]/(x^p)$. One could view this pair as a kind of linearisation of our Hopf algebras, no longer isomorphic to group algebras and group function algebras respectively. Compared to these,  $A_1$ has a modified algebra relation and dually ${A_1^*}$ has a modified coproduct. 

\begin{corollary}\label{Fxt} The coevaluation and the canonical Hopf algebra Fourier transform in terms of $x,L$ take the form
\[\exp=e^{L\tens x}\in {A_1^*}\tens A_1\]
\[ \CF: A_1\to  {A_1^*},\quad \CF(f)=\int e^{L\tens x}f(x),\quad  \CF^{-1}(f)=\int f(L) e^{- L\tens x}.\] 
\end{corollary}
\proof We deduce this as
\[ \exp=\sum_{i=0}^{p-1}(e^L)^i\tens\delta_i(x)= \sum_{m=0}^{p-1}{L^m\over m!}\tens\sum_{i=0}^{p-1}\delta_i(x)i^m=e^{L\tens x}\]
using $x^m=\sum_i\delta_i(x)i^m$. This in turn gives the Fourier transform as stated. In principle, one can also find from (\ref{txpair}) that   $\<x^i,L\>=-\sum_{k=1}^{p-1}{ k^i\over k}=-\sum_{k=1}^{p-1}k^{i-1}=1$ if $i=1$ mod $p-1$ and zero otherwise, to eventually find $\exp$ from this. \endproof

Thus, working with $L$ puts the Fourier transform into a familiar form. We now turn to differentials. We recall that the regular $d=1$ monics are of the form $m=x-\mu$ for $\mu\in \F_p$ (where the corresponding field extension is trivial, so we identify this with $\mu$ in the general construction). The calculus is 1-dimensional with basis $\extd x$ and necessarily descends to $A_1$. However, the classical calculus on $\F_p[x]$ given by $\mu=0$, aside from not being regular, implies $\extd x=\extd x^p=0$ and hence gives the zero calculus on $A_1$. We therefore exclude it in what follows. 

\begin{proposition} For any $\mu\in\F_p^*$, the inherited calculus $\Omega(A_1)$ has $\Omega^i=0$ for $i>1$ and $\Omega^1=A_1\extd x$ with relations 
\[ [\extd x, f]= \mu \extd f.\]
Moreover, $H_{\rm dR}^0(A_1)=\F_p$, $H_{\rm dR}^1(A_1)=\F_p$, spanned by $1$ and $x^{p-1}\extd x$ respectively.
\end{proposition}
\proof The inherited calculus has the form stated, with $\extd f= (\del f)\extd x$ where
\[ \del f={f(x+\mu)-f(x)\over\mu},\quad\forall f\in A_1.\]
We already know $H_{\rm dR}^0$ from Corollary~\ref{H0Ad} but we can also see this directly. If $\del f=0$ then $f(x+\mu)=f(x)$. But $n\mu=\lambda$ mod $p$ has a solution $n$ for all $\lambda$ so by iteration, $f(x+\lambda)=f(x)$ for all $\lambda$. By (\ref{deltaf}), $f(x)$ is determined by its values and we see that these are constant, hence $f$ is a multiple of 1. For $H_{\rm dR}^1$, all 1-forms are closed and if $f\extd x =\extd h$ for some $h(x)$ then $f=\del h=(h(x+\mu)-h(x))/\mu$. Clearly this cannot happen for $f$ of degree $p-1$ since $h$ would need degree $p$ which is not possible. For smaller degree one can iteratively solve to find $h$ by calculations that are the same as for the trivial 1st cohomology of $\F_p[x]$ with its 1-dimensional calculi. The calculus is manifestly inner with $\theta=\mu^{-1}\extd x$. \endproof

Note that calculi on finite sets correspond to directed graphs\cite{Ma:gra} and the above calculi correspond to the Cayley graph on $\Z/p\Z$ generated by singleton sets $\{\mu\}\subset \Z/p\Z$. The directed graph here has edges of the form $i{\buildrel \mu\over\longrightarrow} i+\mu$ corresponding to a  finite difference with step $\mu$ on $\Z/p\Z$. It is easy to see from (\ref{xtF})  by a change of variables in $\CF$ that 
\begin{equation}\label{Foudif} \CF\del f=  \CF(f) \left({t^\mu-1\over \mu}\right),\end{equation}
in keeping with the general features of Fourier transform. 

One can also ask on the dual side about the calculus on ${A_1^*}=\F_p[t]/\<t^p-1\>$. Usually in the abelian case the problem reverts to calculi on the dual group but that is not possible in our case where the order of the group is the characteristic. However, it remains in any characteristic that translation invariant calculi on group algebras are classified by group 1-cocycles\cite{MaTao}. In our case there is a natural choice in which the values of the cocycle are in $\F_p$ with trivial group action. In that case a group cocycle means a group homomorphism from $\Z/p\Z$ to itself, which since $p$ is prime can only be trivial  or the identity. We therefore have a unique 1-dimensional calculus from this point of view, namely
\[  \Omega^1({A_1^*})={A_1^*}v,\quad v=t^{-1}\extd t,\quad \extd t^i=i t\, v,\quad [\extd t,t]=0\]
and $\Omega^2=0$. We see that this is the classical calculus on the algebraic circle $\F_p[t,t^{-1}]$ descended to ${A_1^*}$. Writing $\extd f(t)=(\del f)(t)v$, we have $\del t^m=mt^m$, the degree operator. From (\ref{xtF}) one easily finds
\[  \CF^{-1}\del=x\CF^{-1}\]
so that differentiation on ${A_1^*}$ again becomes multiplication in $A_1$ under Fourier transform. In terms of $L$, we have
\[ v=e^{-L}\extd e^L=e^{-L}\sum_{i=1}^{p-1}{L^{i-1}\over(i-1)!}\extd L=(1+L^{p-1})\extd L.\]
The calculus here descends from the classical calculus on $\F_p[L]$ but $v$ and not $\extd L$ is the basic translation-invariant differential form, because this property depends on the coproduct on $L$ and this was modified from the additive one. Consequently, we have 
\[ \del L=1-L^{p-1},\quad \del L^i=iL^{i-1},\quad \forall i>1\]
for the left-invariant derivative. 

We now ask if there is a quantum Riemannian structure on $A_1$ for the above calculus. For this we must specify the space of 2-forms and the canonical choice here is for $\extd x$ to square to zero, so $\Omega^2=0$.

\begin{proposition} The above $\Omega^1$ on $A_1$ admits a quantum metric  $g$ if and only if $p=2$ and then $g=\extd x\tens\extd x$. This admits only one quantum Levi-Civita connection, given by $\nabla\extd x=0$ and $\sigma(\extd x\tens\extd x)=\extd x\tens\extd x$. 
\end{proposition}
\proof An element of $\Omega^1\tens_{A_1}\Omega^1$ for the above calculus has the form $g=\alpha\extd x\tens\extd x$ for some nonzero element $\alpha\in A_1$. However, a quantum metric to be invertible must also be central and in our case $[g,x]=2\mu g$ which is zero only if $p=2$. Now setting $p=2$, we need $\alpha$ to be invertible in which case $\alpha=1$. Next, we take a general form of connection $\nabla\extd x=a \extd x\tens\extd x$ for $a\in A_1$. If this is a bimodule connection then
\[ \nabla((\extd x)x)=a\extd x\tens(\extd x)x+ \sigma(\extd x\tens\extd x)=\nabla((x+1)\extd x)=\extd x\tens\extd x+(x+1)a\extd x\tens\extd x\]
which requires $\sigma(\extd x\tens \extd x)=(1+a)\extd x\tens\extd x$. This indeed defines a bimodule map as $\extd x\tens\extd x$ is central. All connections are necessarily flat and torsion free due to the choice of $\Omega^2$, so all that remains is metric compatibility. This requires
\begin{align*}\nabla g&=\nabla\extd x\tens\extd x +\sigma(\extd x\tens a \extd x)\tens\extd x=
a\extd x^{\tens3} +\sigma((\extd x)a\tens \extd x)\tens\extd x\\
&=a\extd x^{\tens 3}+(a+\del a)\sigma(\extd x\tens\extd x)\tens\extd x=(a+(a+\del a)(1+a))\extd x^{\tens 3}=0\end{align*}
where the first equality is $\nabla$ applied on the two factors of $g$, with $\sigma$ used to swap the left output of the second instance to the far left. We also used the commutation relations between $\extd x$ and a general element $a$. For the result to vanish, we need $(1+a)\del a=a$,  which is only solved by $a=0$.  
\endproof

This $g=\extd x\tens\extd x$ is translation-invariant as the coefficients in the basis are constant and should be seen as the intrinsic geometry of $A_1$ over $\F_2$, with here only the trivial quantum Levi-Civita connection $\nabla(f\extd x)=\extd f\tens\extd x$ for any $f\in A_1$. The translation-invariant geometry for $A_2$ will be more interesting. 

\subsection{Fourier transform and geometry on $A_2$}  \label{sec:ConCompA2}

Here we consider 
\[ A_2={\F_p[x]/(x^{p^2}-x)}\]
as a Hopf algebra with $x$ primitive. This is isomorphic as a ring to $\F_p^p\times \F_{p^2}^{p(p-1)\over 2}$ and hence is not functions on any finite group. Rather, by Theorem~\ref{cocycle} we know that we can identify
\[ A_2= C_2\tens_\chi A_1,\quad  C_2=\F_p[y]/(y^p+y)\]
for a certain cocycle $\chi$, where $y=g_1(x)$ and $A_1$ is embedded as $\delta_i(x)$. The  structure of $C_2$ is almost that of $A_1$ itself and is exactly $A_1$ if $p=2$. We focus on this simpler case, which is also the only case where the quantum Riemannian geometry is manageable by known methods. In this case $A_2= A_1\tens_\chi A_1\cong\F_2(\Z/2\Z)\tens_\chi\F_2(\Z_2/2\Z)$ where the first copy has generator $y=g_1(x)=x^2+x$ in $A_2$ and the function algebra  description is via Lemma~\ref{A1isom} with the Kronecker $\delta_i$ in the second copy appearing in $A_2$ as $\delta_i(x)$. We write $\bar\delta_i$ for the parallel Kronecker basis of the first copy, embedded as $\bar\delta_0=1+y$ and $\bar\delta_1=y$. The isomorphism in the other direction is provided by the factorisation in $A_2$,
\[  1=\sum_{i,j}\bar\delta_i\delta_j,\quad x=(\bar\delta_0+\bar\delta_1) \delta_1,\quad x^2= \sum_{i\ne j}\bar\delta_i\delta_j,\quad x^3= \sum_{{\rm not}\ i=j=0}\bar\delta_i\delta_j \]
which one can then write as a  tensor product of the factors for the $A_1\tens_\chi A_1$ description. 
The translation-invariant integral for $A_2$ (as for all $A_d$) is required  to have support only on the top degree, which in our case means $\int x^3=1$ and zero on smaller degree monomials. This corresponds under the isomorphism to a tensor product of integrals, so  $\int\bar\delta_i\delta_j=1$ for all $i,j\in \F_2$.  The final ingredient for Fourier transform is a description of the Hopf algebra dual.

\begin{proposition} For $p=2$, the Hopf algebra dual is $A_2^*\cong \F_2[s,t]/(s^2-1,t^2-1)$ as an algebra, with coalgebra and antipode
\[ \Delta s=s\tens st+st\tens s+st\tens st,\quad \Delta t=t\tens t,\quad \eps s=\eps t=1,\quad S s=s,\quad S t=t.\]
Fourier transform  $\CF:A_2\to A_2^*$ is then given by
\[ \CF(1)=1+s+t+st,\quad \CF(x)=(1+s)t,\quad \CF(x^2)=s+t,\quad \CF(x^3)=s+t+st.\]
\end{proposition}
\proof In our new terms, the cocycle and relations in Theorem~\ref{cocycle} are 
\[ \chi(\delta_i\tens\delta_j)=\bar\delta_{i+j+ij},\quad \delta_i^2=\delta_i+\bar\delta_1,\quad\delta_0\delta_1=\bar\delta_1.\]
This description implies that $A_2^*\cong A_1^*\tens^{\chi^*}A_1^*$ as a cocycle coproduct Hopf algebra, where $\chi^*:A_1^*\to A_1^*\tens A_1^*$ is the dualisation of $\chi$. First, as an algebra
\[ A_1^*\tens^{\chi^*}A_1^*\cong\F_2\Z/2\Z\tens\F_2\Z/2\Z\]
by our results in the preceding section. This is the group algebra of $(\Z/2\Z)^2$ and we write it as stated with involutive generators $s,t$. We use the pairing as in (\ref{txpair}) for each copy 
 whereby $\{t^i\}$ and $\{\delta_i\}$ are dual bases and so are $\{s^i\}$, $\{\bar\delta_i\}$. Using this, the cocycle dualises to
\[ \chi^*(1)=1\tens 1,\quad \chi^*(s)=(1+t)\tens (1+t)-1\tens 1\]
after which we use the general Hopf algebra construction \cite[Prop.~6.3.8]{Ma:book}
\[ \Delta(x\tens y)=x\o\tens \chi^*(x\th)\uo y\o\tens x\t\tens \chi^*(x\th)\ut y\t\]
adjoint to (\ref{cocyprod}), where $x\tens y\in A_1^*\tens A_1^*$ is taken on the right hand side with its original tensor product coalgebra and we have written $\chi^*=\chi^*{}\uo\tens\chi^*{}\ut$ (sum understood). This computes for $s,t$ grouplike to the formula stated. The second copy of $A_1^*$ is a sub-Hopf algebra and the first copy is a subalgebra with a cocycle-modified coproduct.  

Once we have the dual Hopf algebra in this form, we have the overall pairing and an integral on $A_2^*$
\[\<\bar\delta_i\delta_j,s^kt^l\>=\delta_{i,k}\delta_{j,l},\quad\int s^it^j=\delta_{i,0}\delta_{j,0}\]
where may check that the latter remains invariant. Then clearly $(\int\tens\int)(\exp)=1$ and
\[ \CF(\bar\delta_i\delta_j)= \int \bar\delta_i\delta_j\bar\delta_m\delta_n\tens s^m t^n=s^it^j  \]
after a short computation using the product in $A_2$ (or the cocycle product from the $A_1\tens_\chi A_1$ point of view). The modified product of $\delta_j\delta_m$ does not affect the answer after integration. In terms of the original description of $A_2$, this comes out as stated. 
\endproof

We also note that the fixed subalgebra  $B_2\subset A_2$  has dimension $3$ according to Proposition~\ref{dimBd}, so  is the Boolean algebra on 3 elements. One has

\begin{proposition} $A_2$ for $p=2$ is reduced and every element obeys $a^4=a$ for all $a\in A_2$.  Moreover, $A_2\cong \F_2 x\oplus B_2$ as a vector space and contains $B_2$ as a subalgebra with orthogonal idempotents $e_1,e_2,e_3$. The Hopf algebra structure of $A_2$ in this form is
\[ e_ie_j=e_i\delta_{ij},\quad\sum_ie_i=1,\quad x^2=x+e_1,\quad e_1x=e_2+x,\quad e_2x=e_2,\quad e_3x=0\]
\[\eps x=\eps e_1=\eps e_2=0,\quad \eps e_3=1,\quad \Delta x=x\tens 1+1\tens x,\quad \Delta e_1=e_1\tens 1+1\tens e_1\]
\[ \Delta e_2=e_2\tens 1+1\tens e_2+e_1\tens x+x\tens e_1,\quad  \Delta e_3=1\tens 1+e_3\tens 1+1\tens e_3+ e_1\tens x+x\tens e_1\]
\[ \eps(e_1)=\eps(e_2)=\eps(x)=0,\quad \eps(e_3)=1.\]
\end{proposition}
\proof By writing $a=\alpha+\beta x+\gamma x^2+\delta x^3$ we see that $a^2=\alpha+\beta x^2+\gamma x+\delta x^3$ and $a^4=a$. This is also clear from the ring structure. The  coefficients here are $0,1$ and in this case $a^n=0$ is not possible for any $n>0$ unless
$a=0$. The boolean elements (meaning $a^2=a$) are of the form $\alpha+\beta(x+x^2)+\delta x^3$ and these form a subalgebra. Here $1, e_1=x^2+x, e_3=x^3+1$ obey $e_1e_3=0$ so with $e_2=1+e_1+e_3=x(x^2+x+1)$ are a complete set of idempotents for this subalgebra. So  $A_2\cong \F_2.x\oplus B_2$.  We easily work out the Hopf algebra structure as stated. The antipode is the identity map.  \endproof

We now turn to the  inherited structure of $\Omega(A_2)$ and its intrinsic translation-invariant geometry. For the calculus, there is in fact only one monic irreducible of degree 2 in $\F_2[x]$ namely $m(x)=x^2+x+1$, so only one such calculus to consider. 

\begin{proposition}\label{OmegaA2} The quotient $\Omega(A_2)$ is 2-dimensional in degree 1 with basis $\extd x,\mu$ and relations
\[ [\extd x,x]=\mu,\quad [\mu,x]=\extd x+\mu.\]
Moreover, the calculus is bicovariant,  inner with $\theta=\extd x+\mu$ and connected with Poincar\'e duality in the sense  
\[ H_{\rm dR}^0(A_2)=\F_2,\quad H_{\rm dR}^1(A_2)=\F_2^2,\quad H_{\rm dR}^2(A_2)=\F_2.\] These are spanned by $1, \{x\extd x, \mu x^2\}$ and $x^3 \extd x\wedge\mu$ respectively. 
\end{proposition}
\proof We work in the $\Omega^1(\F_2[x])=\F_4[x]$ description, reduced to $A_2$, but we write $\extd x=1\in \F_4[x]$ to avoid confusion with $1\in A_2$. Thus $\extd x. x= 1.(x+\mu)=x.1+\mu=\extd x+\mu$. Similarly, $\mu x=(x+\mu)\mu=x\mu+(1+\mu)=\extd x+ (x+1)\mu$ as stated. These are the same basis and relations as the calculation for $\F_2[x]$, just adopted for our quotient algebra. Note also that the calculus necessarily remains inner with $\theta=\mu^{-1}=\extd x+\mu\in \Omega^1$. It necessarily remains bicovariant. If we let $\extd_n$ denote the restriction of the derivative to the $n$-th component of the graded exterior algebra and write $\extd_0 f = \del_1 f  \extd x + \del_2 f  \mu$ for $f \in A_2$, then $\extd_1(f_1  \extd x + f_2  \mu) = (\del_1 f_2 - \del_2 f_1) \extd x \wedge \mu$ for $f_1  \extd x + f_2  \mu \in \Omega^1$.  We already know $H_{\rm dR}^0$ by Corollary~\ref{H0Ad} but it is also easy to verify directly.  Brute-force calculation shows that $\text{Im}(\extd_0)$ is spanned over $\F_2$ by $\{\extd x, \mu, x^2\extd x + x \mu \}$, $\text{ker}(\extd_1)$ is spanned by $\{\extd x,\mu,x\extd x, x^2\extd x+x \mu, x^2 \mu\}$, and finally that $\text{Im}(\extd_1)$ is spanned by $\{1,x,x^2\}\extd x\wedge\mu$.   The dimensions and bases of the cohomologies follow. Note that over $\F_2$ the exterior algebra is both commutative and anticommutative and symmetric combinations of the basic 1-forms are in the kernel of $\wedge$.
 \endproof

By contrast, this cohomology does not hold for the universal calculus on $A_2$ which is necessarily acyclic and hence cannot obey Poincar\'e duality, and has weaker relations 
\[ [\extd x,x]=\mu,\quad [\mu,x]=\theta,\quad [\theta,x]=\extd x\]
where $\theta$ is an independent 1-form. The  $\Omega(A_2)$ in Proposition~\ref{OmegaA2} is the quotient of this by a further relation $\theta=\extd x+\mu$ which respects the coaction so that the result remains bicovariant.   

We now turn to the quantum Riemannian geometry with this inherited 2-dimensional calculus. Note that in noncommutative geometry a metric, when it exists, need not admit a `Levi-Civita' connection (in the sense of torsion free and metric compatible) and if it does, the connection need not be unique.  In the Hopf algebra case it is natural to consider left-invariant metrics, i.e. ones that are constant in the basic 1-forms, in our case $\extd x,\mu$. 

\begin{proposition}\label{F2levi} There are three left-invariant quantum metrics  $g\in \Omega^1\tens_{A_2}\Omega^1$ namely 
\[ g=\alpha(\mu\tens\mu+\theta\tens\theta)+\beta(\extd x\tens\extd x+\theta\tens\theta)\]
where $\alpha,\beta\in \F_2$ (at least one of them nonzero), each with precisely two invariant torsion free metric compatible bimodule connections, namely $\nabla\extd x=\nabla\mu=0$, $\sigma={\rm flip}$ on the generators and
 \[ \nabla\extd x=\alpha\extd x\tens\extd x+\beta(\extd x\tens\mu+\mu\tens\extd x)+\alpha\beta\mu\tens\mu\]
 \[ \nabla\mu=\beta\mu\tens\mu+\alpha(\extd x\tens\mu+\mu\tens\extd x)+\alpha\beta\extd x\tens\extd x\]
 \[ \sigma(\extd x\tens\extd x)=\alpha\extd x\tens\extd x+\beta \theta\tens\theta+ \alpha\beta(\extd x\tens\mu+\mu\tens\extd x)\]
 \[ \sigma(\mu\tens\mu)=\beta\mu\tens\mu+\alpha \theta\tens\theta+ \alpha\beta(\extd x\tens\mu+\mu\tens\extd x)\]
 \[ \sigma(\extd x\tens\mu)=\alpha\theta\tens\extd x+\beta\mu\tens\theta+\alpha\beta(\mu\tens\extd x+\theta\tens\theta)\]
 \[ \sigma(\mu\tens\extd x)=\alpha\extd x\tens \theta+\beta\theta\tens \mu+\alpha\beta(\extd x\tens \mu+\theta\tens\theta).\]
Moreover, the connections in both cases are flat. 
\end{proposition} 
 \proof We let $f=\extd x\tens \mu + \mu\tens\extd x$ and $h=\extd x\tens\extd x+\mu\tens \mu$ and compute
\[ [\extd x\tens\extd x, x]=f=[\mu\tens\mu,x],\quad [\extd x\tens\mu,x]=\extd x\tens\mu +h,\quad [\mu\tens\extd x,x]=\mu\tens\extd x+h\]
from which it follows that central combinations must be of the form 
\[ g=\alpha \extd x\tens\extd x+\beta\mu\tens\mu +(\alpha+\beta) f\]
which can also be written as stated. Here $\alpha,\beta$ could be functions. We now focus on the constant case. Writing in the $\extd x,\mu$ basis, invertibility then needs $\alpha\beta+\alpha+\beta=1$ which is all cases except $\alpha=\beta=0$. 
 
 Next, we look for bimodule connections. The direct approach is not practical and we use a result in \cite{Ma:gra} that when the calculus is inner, as it is here with $\theta=\extd x+\mu$, bimodule connections are of the form $\nabla \omega=\theta\tens\omega-\sigma(\omega\tens\theta)+\tilde\alpha\omega$ for bimodule maps $\sigma,\tilde\alpha$, and are torsion free if and only if 
\begin{equation}\label{torinner} \wedge\tilde\alpha=0,\quad \wedge\sigma=-\wedge\end{equation}
and metric compatible if and only if
\begin{equation} \label{ginner} \theta\tens g +(\tilde\alpha\tens\id) g+\sigma_{12}(\id\tens(\tilde\alpha-\sigma_\theta)) g=0.\end{equation}
 Thus, if we suppose a map 
 \[ \tilde\alpha(\extd x)=a\extd x\tens\extd x+ b\theta\tens\theta+ c f    \] 
 then 
 \[ \tilde\alpha(\mu)=\tilde\alpha([\extd x,x])=[\tilde\alpha(\extd x),x]=(a+b+c)f\]
 using $[f,x]=f$ and the above. Then
 \[ a\extd x\tens\extd x+b\mu\tens\mu+(a+b)f=\tilde\alpha(\extd x+\mu)=\tilde\alpha([\mu,x])=[\tilde\alpha(\mu),x]=(a+b+c)f\]
 requires $a,b,c=0$. Hence there are no non-zero module maps $\tilde\alpha$ with the required property in (\ref{torinner}). We therefore drop the bimodule map $\tilde\alpha$. Similarly, let
 \[ \sigma(\extd x\tens\extd x)=a\extd x\tens \extd x+b\mu\tens\mu + c f\]
 \[ \sigma(\mu\tens\mu)=A\extd x\tens \extd x+B\mu\tens\mu + C f\]
 \[\sigma(\extd x\tens \mu)=a' \extd x\tens \extd x+ b' \mu\tens\mu+c' f+ \mu\tens\extd x\]
 as dictated by $\wedge\sigma=-\wedge$. Then
 \[ \sigma(f)=\sigma([\extd x\tens\extd x,x])=[\sigma(\extd x\tens\extd x),x]=(a+b+c)f=(A+B+C)f\]
  (by $f=[\mu\tens\mu,x]$ for the second version) so that 
 \[ a+b+c=A+B+C.\]
Similarly $\sigma([\extd x\tens\mu,x])=[\sigma(\extd x\tens\mu),x]$ gives us two further equations
 \[ a'=1+a+A,\quad b'=1+b+B.\]
 This leaves us parameters  $a,b,c,A,B,c'$ for $\sigma$ with the required symmetry. Then writing $\sigma_\theta=\sigma(\id \tens\theta)$ we have
 \[\sigma_\theta(\extd x)=(1+A)\extd x\tens\extd x+(1+B)\mu\tens\mu+(c+c')f+\mu\tens\extd x,\]
\[ \sigma_\theta(\mu)=(1+a)\extd x\tens\extd x+(1+b)\mu\tens\mu+(A+B+c')f+\mu\tens\extd x\]
and hence
\[ \nabla\extd x=A\extd x\tens\extd x+(1+B)\mu\tens\mu+(c+c')f,\]
\[ \nabla\mu=(1+a)\extd x\tens\extd x+b \mu\tens\mu+(1+A+B+c')f.\]
Up to this point we have been fairly general but we now assume the connection is translation-invariant which amounts to our functions being constants. Then 
\begin{eqnarray*}(\id\tens\sigma_\theta)g&=&\extd x^{\tens 2}\left((\alpha(a+A)+\beta(1+a))\extd x+(\alpha(A+B+c)+\beta(A+B+c'))\mu\right)\\
 &+&\mu^{\tens 2}\left((\alpha(1+c+c')+\beta(A+B+c))\extd x+(\alpha(1+B)+\beta(b+B))\mu\right)\\
 & +&\extd x\tens\mu\left((\alpha(A+B+c)+\beta(1+A+B+c'))\extd x+(\alpha(b+B)+\beta(1+b))\mu\right)\\
 & +&\mu\tens\extd x\left((\alpha(1+A)+\beta(a+A))\extd x+((\alpha(c+c')+\beta(A+B+c))\mu\right).\end{eqnarray*}
Applying $\sigma\tens\id$ and equating to $\theta\tens g$ so as to solve the  metric compatibility equation (\ref{ginner}) we obtain a system of quadratic equations for our 6 parameters. Over $\F_2$, we try all 64 parameter values for each of the three non-zero cases of $\alpha,\beta$, finding two solutions in each case. These are the unique nontrivial connections stated and one common connection which is zero on the basic forms and for which $\sigma$ flips the generators as is the case classically. One may then verify metric compatibility directly as a check. That all four connections have zero curvature is obvious for the trivial one and a calculation for the other case. For example
\begin{align*}R_\nabla \extd x&=(\extd\tens\id-\id\wedge\nabla)\nabla\extd x=(\alpha\extd x+\beta\mu)\wedge\nabla\extd x+(\alpha\beta\mu+\beta\extd x)\wedge\nabla\mu=0\end{align*}
where we used the solution for $\nabla\extd x$ and the $\extd\tens\id$ does not contribute as all the coefficients are constant. Using $\nabla\extd x$ and $\nabla\mu$ and that only $\mu\wedge\extd x=\extd x\wedge\mu$ products are non-zero and collecting $\extd x\wedge\mu\tens\extd x$ and $\extd x\wedge\mu\tens \mu$ terms, we obtain zero. Similarly for $R_\nabla\mu=0$. 
 \endproof
 
 The trivial connection here can still be nonzero since $\nabla(a\extd x+b\mu)=\extd a\tens\extd x+\extd b\tens\mu$ for all $a,b\in A_2$, and corresponds geometrically to what we might expect on an affine line. The other connection in each case is more unexpected and it is remarkable that for each metric we find a unique other one. The existence of such a second `nonclassical' quantum Levi-Civita connection was also a feature in the concrete model in \cite{BegMa4}. The general case of nonconstant $\alpha,\beta$ and non-constant connection coefficients in Proposition~\ref{F2levi} is much harder but can in principle be analysed in the same way with additional $\extd\alpha,\extd\beta$ terms entering in the equations for the connection.

\end{document}